\documentclass[reqno]{alea2}
\usepackage{natbib}
\usepackage{fancyhdr}
\usepackage{graphicx}

\usepackage{amssymb}
\usepackage{amsmath}

\renewcommand{\cite}{\citet}

\makeatletter \@addtoreset{equation}{section} \makeatother

\renewcommand\theequation{\thesection.\arabic{equation}}
\renewcommand\thefigure{\thesection}

\renewcommand\thetable{\thesection.\@arabic\c@table}

\newtheorem{theorem}{Theorem}[section]

\newtheorem{proposition}[theorem]{Proposition}

\newtheorem{definition}[theorem]{Definition}

\newcommand{\mc}[1]{{\mathcal #1}}

\usepackage{latexsym,amssymb,amsmath,epsfig}
\usepackage{a4wide}
\usepackage{graphicx}

\def\RR{\mathbb{R}}
\def\NN{\mathbb{N}}
\def\ZZ{\mathbb{Z}}

\def\EE{\mathbb{E}}

\renewenvironment{proof}{\vskip 2mm\noindent {\it Proof:}}
                    {\hfill $\square$ \vskip 2mm \noindent}

\begin{document}
\date{5/20/06, revised}

\author[1]{David Coupier}
\author[2]{Paul Doukhan}
\author[3]{Bernard Ycart}

\address[1]{MAP5, UMR CNRS 8145, Universit\'e Ren\'e Descartes, Paris}
\address[2]{LS-CREST, URA CNRS 220, Paris}
\address[3]{LMC-IMAG, UMR CNRS 5523, Grenoble}

\email[1]{coupier@math-info.univ-paris5.fr}
\email[2]{doukhan@ensae.fr}
\email[3]{ycart@imag.fr}

\urladdr{http://www-lmc.imag.fr/lmc-sms/Bernard.Ycart/publis/zufields.pdf}

\title[Zero-one laws for random fields]{Zero-one laws for
binary random fields} 

\begin{abstract}
A set of binary random variables indexed by a lattice torus is
considered. Under a mixing hypothesis,
the probability of any proposition belonging to the first order
logic of colored graphs tends to $0$ or $1$,  as the size
of the lattice tends to infinity.
For the particular case of the Ising model with bounded pair
potential and surface potential tending to $-\infty$,
the threshold functions of local propositions are computed,
and sufficient conditions for the zero-one law are given.
\end{abstract}

\maketitle
\section{Introduction}
\label{intro}
There exist essentially three types of zero-one laws. The most ancient
(Kolmogorov, Hewitt-Savage) concern asymptotic events defined on a
sequence of independent random variables (see \cite{Feller2}
p.~124). More recently, results describing how the probability of an
increasing subset of $\{0,1\}^n$ goes from almost $0$ to almost $1$
have been developed (see for instance \cite{Talagrand94a}). Zero-one laws
in logic state that all propositions of a given logic have a
probability tending to $0$ or $1$, as the size of the domain to which
they apply tends to infinity. They have a long history: see
\cite{Compton} or \cite{Winkler} for reviews,
and \cite{Spencer01} on logic and random graphs.
We shall use \cite{EbbinghausFlum} as a basic
reference on first order logic.

The first zero-one law in logic was proved
independently by \cite{glebskiietal} and
\cite{Fagin}. It states that the probability of any first order
property on a finite domain without constraints tends to $0$ or $1$ as
the size of the domain tends to infinity, under the
uniform distribution on the set of structures.
It was soon recognized that it also holds
when the set of structures is endowed with a product of
Bernoulli distributions with fixed parameter $p$. As an elementary
example, consider the model made of a single unary
predicate, applied to the domain $V_n = \{0,\ldots,n\!-\!1\}$
(see \cite{Spencer01}, Section 0.2). Any structure with $k$ asserted
and $n\!-\!k$ negated facts has probability $p^k(1-p)^{n-k}$. The
Glebskii et al.--Fagin theorem applies, but the expressive power
of the language is very poor. Let us add a binary predicate
$R$. The elements of $V_n$ are interpreted as vertices of a graph, and
any fact $Rxy$ will be seen as an existing edge between vertices $x$
and $y$. The unary predicate, denoted by $C$, will  be
interpreted as a coloring of
vertices: the fact $Cx$ means that $x$ is a black vertex and $\neg Cx$
that it is white. Among all possible $2^{n+n^2}$ structures, 
we shall retain only those $2^n$ for which the graph is
cyclic:
$$
Rxy \Longleftrightarrow y=x\pm 1 \mod n\;.
$$
The expressive power of the language has now increased, since it allows
statements about the geometrical arrangement of colored vertices. For
instance the first order logic now contains the sentence ``there
exist $3$ black vertices in a row''.
 This model is essentially the ``random circular unary predicate''
model of \cite{Spencer93}, Section 2.
One can see it as a fixed graph structure
(the $n$-cycle) with randomly colored vertices. The random colors
are i.i.d. random variables, Bernoulli with parameter $p$. The
Glebskii et al.--Fagin theorem does not apply anymore, since the set of
structures has been restricted. However, it is easy to see that the
zero-one law still holds. But until now only {\it  independent}
random colors have been considered. The aim of this article is to
study situations under which the zero-one law holds even though
the colors of vertices are dependent random variables.

Let us generalize the previous model by considering a lattice graph
in dimension $d$, with periodic boundary conditions (lattice torus).
The vertex set is now $V_n = \{0,\ldots n\!-\!1\}^d$. The edge set, denoted
by $E_n$, will be specified by defining the set of neighbors ${\mc V}(x)$
of a given vertex $x$:
\begin{equation}
\label{def:voisinage}
{\mc V}(x) = \{ y\neq x\in V_n\,,\; \|y-x\|_p \leq \rho\}\;,
\end{equation}
where the substraction is taken componentwise modulo $n$,
$\|\cdot\|_p$ stands for the $L_p$ norm in $\RR^d$
($1\leq p\leq\infty$), and $\rho$ is a fixed parameter.
For instance, the square lattice is obtained for $p=\rho=1$. Replacing
the $L_1$ norm by the $L_\infty$ norm adds the diagonals.

A lexical bridge between probability and logic is needed here.
A {\it binary random field} is a set of random
variables indexed by $V_n$, with values in $\{-1,+1\}$
(see \cite{Rosenblatt} as a general reference). Obviously the
choice $\{-1,+1\}$ for the values of the field is arbitrary:
it could be replaced by $\{0,1\}$ or $\{\mbox{white},\mbox{black}\}$.
A binary random field can also be
viewed as a random mapping from $V_n$ to $\{-1,+1\}$, which is
usually called a {\it configuration} (a spin configuration for 
the Ising model).
The set of configurations will be denoted by $\Xi_n$.
In the logical interpretation, the graph structure $(V_n,E_n)$
can be described using binary predicates (see Section \ref{logic} for
details). The values of a configuration can be viewed as a coloring of
the graph, and described by a unary predicate $C$. Let $\eta$ be a
configuration. The assertion  $\eta(x)=+1$ will be identified to the
fact $Cx$, and $\eta(x)=-1$ to its negation $\neg Cx$. For convenience,
we shall also retain the coloring interpretation:
$$
\begin{array}{ccccl}
\eta(x)=+1&\Longleftrightarrow &Cx&
\Longleftrightarrow& \mbox{vertex $x$ is black}
\\
\eta(x)=-1&\Longleftrightarrow &\neg Cx&
\Longleftrightarrow& \mbox{vertex $x$ is white}
\end{array}
$$
Thus the set of configurations $\Xi_n$ will be identified to the set of
those logical structures on the domain $V_n$ for which facts relative
to binary predicates are fixed and compatible with the chosen lattice
graph. Any closed logical formula (also called sentence) $A$,
expressed in terms of the predicates,
determines a subset $A_n$ of $\Xi_n$ that satisfies it.
Given a probability distribution $\mu_n$
on $\Xi_n$, the zero-one law holds if $\mu_n(A_n)$ tends to $0$ or $1$,
for each sentence $A$ belonging to the first order logic.
Two hypotheses will be used. The first one
says that the probability of a local
configuration should remain bounded away from $0$. The
second one is a classical mixing condition: the
covariance between functions depending on two disjoint subsets of
vertices tends to zero as the distance between the two subsets
tends to infinity.
Under these hypotheses, Theorem \ref{th:zuweak} shows that
the zero-one law holds. The proof is based on Gaifman's
theorem (\cite{Gaifman82} and \cite{EbbinghausFlum} p.~31), which says
that first-order sentences are essentially local. The
zero-one law for first order sentences can be reduced to computing the
probability of local patterns (see Section 2 for precise definitions
and statements). It will be proved that any local configuration
appears with probability tending to $1$, an analogue of the ``typing
monkey'' paradox.

Quoting \cite{Spencer93}, ``To people who
work on Random Graphs the cases $p$ constant are only a small and
relatively uninteresting part of the theory''. In the circular
model above, when $p=p(n)$ depends on the number of vertices,
threshold phenomenons similar to those of random graphs occur. For
instance, the probability of having $k$ consecutive black vertices
tends to $0$ if $p(n)\ll n^{-\frac{1}{k}}$, it tends to $1$ if
$p(n)\gg n^{-\frac{1}{k}}$. An analogous result
holds for random colorings of a lattice graph in dimension $d$: the
probability of having $k$ neighboring black vertices tends to $0$ or
$1$ according to whether $p(n)$ is small or large compared to the
threshold function $p(n)=n^{-\frac{d}{k}}$ (see
\cite{dcadby1} for an interpretation in the context of random images).
In order to extend it to weakly dependent random fields, one needs to
choose a parametric family of random field distributions. One of the simplest
and most widely studied is the Ising model (see
e.g. \cite{Georgii,KindermannSnell}).
\begin{definition}
\label{def:ising}
Let $G=(V,E)$ be an undirected graph structure with finite vertex set
$V$ and edge set $E$. Let $a$ and $b$ be two reals. The
{\rm Ising model} with parameters $a$ and $b$ is the probability
measure $\mu_{a,b}$ on $X=\{-1,+1\}^V$ defined by:
$$
\forall \eta\in X\;,\quad
\mu_{a,b} (\eta) = \frac{1}{Z_{a,b}}
\exp\left(a\sum_{x\in V} \eta(x)
+b\sum_{\{x,y\}\in E} \eta(x)\eta(y)\right)\;,
$$
where the normalizing constant $Z_{a,b}$ is such that
$\sum_{\eta\in E} \mu_{a,b}(\eta) = 1$.
\end{definition}
In the classical presentation of statistical physics, the elements of
$X$ are spin configurations; the definition of $\mu_{a,b}$ involves a
temperature parameter which is not relevant
here and will be omitted. The parameters $a$ and $b$ are respectively
the ``surface'' and ``pair'' potentials. For $a>0$ (respectively
$a<0$), the measure $\mu_{a,b}$ gives a higher weight to those
configurations with a large number of $+1$'s
(respectively $-1$'s). For $b>0$, the measure
$\mu_{a,b}$ tends to favor groups of neighboring vertices with the
same spin, whereas for $b<0$, any $+1$ will be more likely
surrounded by $-1$'s. For $b=0$, $\mu_{a,0}$ is a product
measure: the spins of all vertices are mutually
independent, $+1$ with probability $p=e^{a}/(e^a+e^{-a})$ and $-1$
with probability $1-p$.
 Observe that the model remains unchanged by swapping $+1$
and $-1$ and replacing $a$ by $-a$. In order to keep a certain
coherence with random graphs, we chose to study negative values
of $a$ (corresponding to small values of $p$).

Let us now consider the Ising model $\mu_{a,b}$ on the lattice graph
$(V_n,E_n)$ defined above. The potentials $a=a(n)$ and $b=b(n)$
depend on the size of the lattice. The interaction potential $b(n)$
controls the degree of local dependence: we will assume that it
remains bounded. Theorem \ref{th:zuising} expresses the threshold
functions of local properties in terms  of
$e^{a(n)}n^{\frac{d}{2k}}$: if this quantity
tends to $0$ or $+\infty$ for all $k$,
then the zero-one law holds.

The paper is organized as follows. In Section \ref{logic} the logical
setting is described, and the zero-one law for first order logic is
reduced to the study of basic local and pattern sentences (Proposition
\ref{prop:reduction}). Weakly dependent random fields will be
introduced in Section \ref{fields}. It will be proved that the
probability of any pattern sentence tends to $1$,
which implies the zero-one law (Theorem
\ref{th:zuweak}). Section \ref{ising} is devoted to the Ising model
with varying potentials $a$ and $b$. A description of
threshold functions for basic local sentences will be given, and
the zero-one law will be deduced.
\section{First order logic}
\label{logic}
We follow the notations and definitions in Chapter $0$ of
\cite{EbbinghausFlum} for the syntax and semantics of first order
logic. The vocabulary contains one unary predicate, denoted by $C$, and
some binary predicates. They apply to the domain
$V_n=\{0,\ldots n\!-\!1\}^d$. Once the domain and the vocabulary are
fixed, the structures are particular models of the predicates,
applied to variables in the domain. To any
structure, a graph is naturally associated (\cite{EbbinghausFlum} p.~26),
connecting those pairs of elements $\{x,y\}$ which are such that
$Rxy$ or $Ryx$ are satisfied, where $R$ is any of the binary predicates.
In each of our structures, we decide to fix the list of facts relative
to binary predicates, so that the associated graph is the lattice
torus defined by (\ref{def:voisinage}).
It will be denoted by $G_n=(V_n,E_n)$.
As usual, the graph distance $dist$ is defined as the minimal
length of a path between two vertices.
The ball of center $x$ and radius $r$ will be denoted by $B(x,r)$.
$$
B(x,r) = \{\, y\in V_n\,;\; dist(x,y)\leq r\,\}\;.
$$
In order to avoid particular cases due to self-overlapping balls, we
will always assume that $n>2\rho r$. If $n$ and $n'$ are both larger
than $2\rho r$, the balls $B(x,r)$ in $G_n$ and $G_{n'}$ are
isomorphic. Two properties of the balls $B(x,r)$ will be crucial in
what follows. The first one is that two balls with same radius are
translated of each other:
$$
B(x+y,r)=y+B(x,r)\;.
$$
(Recall that all operations on vertices are modulo $n$.)
The second one is that for $n>2\rho r$, the cardinality of $B(x,r)$
depends on $r$, but neither on $x$ nor on $n$;
it will be denoted by $\beta(r)$. Observe that whatever the choices of
$p$ and $\rho$ in (\ref{def:voisinage}) the ball $B(x,r)$ is included
in the sub-lattice $[x-\rho r,x+\rho r]^d\cap V_n$, and thus
$$
\beta(r) \leq (2\rho r+1)^d\;.
$$

In our setting, two structures may differ
only by their unary facts. There is a
natural bijection between the set of structures and
$\{-1,+1\}^{V_n}$: to a structure $S$ corresponds the configuration
$\eta$ defined by:
$$
(S\models Cx) \Longleftrightarrow (\eta(x)=+1)\quad;\quad
(S\models \neg Cx) \Longleftrightarrow (\eta(x)=-1)\;.
$$

Formulas such as $Cx$, $Rxy$\ldots~ are called {\it atoms}. The
{\it first-order logic} (\cite{EbbinghausFlum} p.~5)
is the set of all formulas obtained by recursively combining
first-order formulas, starting with atoms.
\begin{definition}
The set ${\mc L}_1$ of first-order formulas is defined by:
\begin{enumerate}
\item
All atoms belong to ${\mc L}_1$.
\item
If $A$ and $B$ are first-order formulas, then
$(\neg A)$, $(\forall x Ax)$ and $(A\wedge B)$ also belong to
${\mc L}_1$.
\end{enumerate}
\end{definition}
Observe that the definition of formulas through variables and atoms
rules out any reference to a particular vertex. For instance the
sentence ``the origin is black'' does not belong to the first order
logic. Moreover,
we are interested in formulas for which it can be decided if they are
true or false for any given configuration, i.e. for which all variables are
quantified. They are called closed formulas, or sentences.
Such a sentence $A$ defines a subset $A_n$ of $\Xi_n$: that of
all configurations $\eta$ that {\it satisfy} $A$ ($\eta\models A$).
$$
A_n = \{\eta\in \Xi_n\,;\;\eta\models A\}\;.
$$
If $\mu_n$ denotes a probability measure on $\Xi_n$, then any sentence $A$
has probability $\mu_n(A_n)$ to be satisfied.
$$
\mu_n(\eta\models A) = \mu_n(A_n)\;.
$$
Until now, we have not precised which binary predicates belong to the
vocabulary. As an example, consider the square torus lattice in
dimension $2$. Each vertex $x=(i,j)$ has $4$ neighbors
$(i\pm 1,j), (i,j\pm 1)$ (additions are always taken modulo $n$). This graph
structure can be described by a single binary predicate $N$ (neighbor):
$$
Nxy \Longleftrightarrow y=x\pm (0,1) \mbox{ or } y=x\pm (1,0)\;.
$$
One can also consider $2$ binary
predicates $R$ (right) and $U$ (up).
$$
Rxy \Longleftrightarrow y=x+(1,0) \mbox{ and }
Uxy \Longleftrightarrow y=x+(0,1)\;.
$$
The underlying graph is the same but the language is more expressive.
The sentence ``there exist three neighboring black vertices'' belongs
to both first order logics, but ``there exist three horizontally
adjacent black vertices'' only belongs to the second one. It turns out
that the zero-one law depends only on the graph, and not on the
binary predicates. This is essentially a consequence of Gaifman's theorem
(\cite{EbbinghausFlum} p.~31), which states
that every first-order sentence is equivalent to a boolean combination of
{\it basic local sentences}.
\begin{definition}
\label{def:basiclocal}
A {\it basic local sentence} has the form:
\begin{equation}
\label{deflocalebasique}
\exists x_1\ldots\exists x_m\;
\left(\bigwedge_{1\leq i<j\leq m} dist(x_i,x_j)>2r\right)
\;\wedge\;
\left(\bigwedge_{1\leq i \leq m} \psi_i(x_i) \right)\;,
\end{equation}
where:
\begin{itemize}
\item[$\bullet$] $m$ and $r$ are fixed nonnegative integers,
\item[$\bullet$] for all $i=1,\ldots,m$,
$\psi_i(x) \in {\mc L}_1$ is a formula for which only variable $x$ is
free (not bound by a quantifier), and the other variables
all belong to the ball $B(x,r)$.
\end{itemize}
\end{definition}
For fixed radius $r$,
consider a {\it complete description}
$D(0)$ of the ball $B(0,r)$, i.e. a first-order sentence for which
all statements concerning vertices at distance at most $r$ of $0$ are
either asserted or negated. It describes a certain local
configuration of the vertices in $B(0,r)$.
Denote this configuration by $\zeta_{D,0}$.
For each $x\in V_n$, denote by $\zeta_{D,x}$ the translation of
$\zeta_{D,0}$ onto the ball $B(x,r)$:
$$
\forall y\in\ZZ^d\,,\; dist(0,y)<r \Longrightarrow
\zeta_{D,x}(x+y) = \zeta_{D,0}(y)\;.
$$
The complete description of the ball $B(x,r)$ that describes the local
configuration $\zeta_{D,x}$ will be denoted by $D(x)$.
\begin{definition}
A {\it single pattern sentence} has the form:
\begin{equation}
\label{def:singlepattern}
\exists x\; D(x)\;,
\end{equation}
where $D(0)$ is a complete description of the ball $B(0,r)$ for a
fixed $r$.
\end{definition}
Examples of single pattern sentences are:
\begin{enumerate}
\item
``there exists a black vertex'',
\item
``there exists a ball of radius $2$ with a black vertex at the center,
  all other vertices in the ball being white''.
\end{enumerate}
\begin{definition}
A {\it pattern sentence} has the form:
\label{def:pattern}
\begin{equation}
\label{defpattern}
\exists x_1\ldots\exists x_m\;
\left(\bigwedge_{1\leq i<j\leq m} dist(x_i,x_j)>2r\right)
\;\wedge\;
\left(\bigwedge_{1\leq i \leq m} D_i(x_i) \right)\;,
\end{equation}
where:
\begin{itemize}
\item[$\bullet$] $m$ and $r$ are fixed nonnegative integers,
\item[$\bullet$] for all $i=1,\ldots,m$,
$D_i(0)$ is a complete description of the ball $B(0,r)$.
\end{itemize}
\end{definition}
Examples of pattern sentences are:
\begin{enumerate}
\item
``there exist $3$ black vertices'',
\item
``there exist $3$ non overlapping white balls of radius $2$, one of
  them with a black vertex at the center''.
\end{enumerate}
Obviously, pattern sentences are particular cases of basic
local sentences. Proposition \ref{prop:reduction} below reduces the
proof of zero-one laws for random configurations to pattern sentences.
The set of predicates being fixed, let $\mu_n$ be any probability
distribution on the set of configurations $\Xi_n$.
\begin{proposition}
\label{prop:reduction}
Consider the following three assertions.
\begin{itemize}
\item[(i)] The probability of any {\rm pattern} sentence tends to $0$ or $1$.
\item[(ii)] The probability of any {\rm basic local} sentence
tends to $0$ or $1$.
\item[(iii)] The probability of any {\rm first order} sentence
tends to $0$ or $1$.
\end{itemize}
Then {\it (i)} implies {\it (ii)} and {\it (ii)} implies {\it (iii)}.
\end{proposition}
\begin{proof}
Observe first that if the probabilities of sentences $A$ and $B$
tend to $0$ or $1$, then so do the probabilities of $\neg A$ and
$A\wedge B$. This follows from elementary properties of
probabilities. As a consequence, if the probability of $A$ tends
to $0$ or $1$ for any $A$ in a given family, this remains true
for any finite Boolean combination of sentences in that family.
Thus Gaifman's theorem yields that {\it (ii)} implies {\it (iii)}.
We shall prove now that every basic local sentence is either
unsatisfiable or a finite Boolean combination of pattern
sentences. Indeed, consider a formula $\psi(x)$ for which only
variable $x$ is free, and the other variables all belong to the
ball $B(x,r)$. Either it is not satisfiable, or there exists a
finite set of local configurations (at most $2^{\beta(r)}$) which
satisfy it. To each of those configurations corresponds a complete
description $D(x)$ which implies $\psi(x)$. So $\psi(x)$ is
equivalent to the disjunction of these $D(x)$'s:
\begin{equation}
\label{decompositionlocalcomplete}
\psi(x) \leftrightarrow \bigvee_{D(x)\rightarrow \psi(x)} D(x)\;.
\end{equation}
In formula (\ref{deflocalebasique}), one can replace each $\psi_i(x_i)$
by a disjunction of complete descriptions. Rearranging terms, one gets
that the basic local sentence (\ref{deflocalebasique}) is itself a
finite disjunction of pattern sentences.
\end{proof}
In the case of weakly dependent random fields, to be treated in
the next section, we will prove that the probability of any
pattern sentence tends to $1$. Proposition \ref{prop:reduction2}
shows that it is enough to examine single pattern sentences. As
Proposition \ref{prop:reduction}, it holds for any probability
distribution $\mu_n$ defined on the set of configurations $\Xi_n$. 
\begin{proposition}
\label{prop:reduction2}
If the probability of any single pattern sentence tends to $1$, then
the probability of any pattern sentence tends to $1$.
\end{proposition}
\begin{proof}
Consider the pattern sentence defined by (\ref{defpattern}). We show
that its probability tends to $1$ by contructing a single
pattern sentence that implies it. Consider
the ball $B(0,R)$, with $R=m(\rho r+1)$. This ball contains $m$
{\it  disjoint} balls of radius $r$.
Denote by $x_1(0),\ldots,x_m(0)$
their centers.

The ball $B(0,R)$ is partitioned into $m+1$ sets: the
$m$ balls $B(x_1(0),r),\ldots,B(x_m(0),r)$ and $\overline{B}(0)$:
$$
\overline{B}(0)=B(0,R)\setminus\bigcup_{1\leq i\leq m} B(x_i(0),r)\;.
$$
Define a complete
description $\tilde{D}(0)$ on $B(0,R)$ as follows:
\begin{enumerate}
\item
For $i=1,\ldots,m$ the restriction of $\tilde{D}(0)$ to $B(x_i(0),r)$
is $D_i(x_i(0))$.
\item
For all $y\in \overline{B}(0)$, $\tilde{D}(0) \models Cy$ (all
vertices outside the balls $B(x_i(0),r)$ are black).
\end{enumerate}
Clearly, $\exists x\;\tilde{D}(x)$ implies (\ref{defpattern}).
By hypothesis, the probability of $\exists x\;\tilde{D}(x)$
tends to one, hence the result.
\end{proof}
\section{Weakly dependent random fields}
\label{fields}
The notations are those of the previous sections. For each positive
integer $n$, we consider a probability measure $\mu_n$ on the set of
configurations $\Xi_n$ on the vertices of a lattice graph
$G_n=(V_n,E_n)$ in dimension $d$.
Theorem \ref{th:zuweak} below gives conditions
under which the probability of any single pattern sentence tends to
$1$, which implies the zero-one law by Propositions
\ref{prop:reduction} and \ref{prop:reduction2}.
The first condition says that the probability of any
local configuration remains bounded away from $0$. The second
one bounds the covariance of two distant configurations. Both
conditions, together with the conclusion of the theorem can be
expressed in terms of probabilities of events such as
``the restriction of the (random) configuration to a given subset
coincides with a given local configuration''.
If $\eta$ is a (global) configuration on $V_n$,
its restriction to a subset $B$ of $V_n$
will be denoted by $\eta_{B}$.
The cardinality of $B$ is denoted by $|B|$. If $B$ and $C$ are two
vertex subsets, their distance is defined as usual by:
$$
dist(B,C) = \min\{\,dist(x,y)\,,\; x\in B\,,\;y\in C\,\}\;.
$$
We call {\it local event relative to $B$} an event which depends only
on the values of $\eta(x)$ for $x\in B$.
\begin{theorem}
\label{th:zuweak}
For each positive integer $n$, let $\mu_n$ be a probability measure on $\Xi_n$.
Assume that the following hypotheses hold.
\begin{enumerate}
\item
Let $r$ be a fixed integer. For any local configuration
$\zeta_{D,0}$ on the ball $B(0,r)$, let $\zeta_{D,x}$ be defined by
$\zeta_{D,x}(x+y)=\zeta_{D,0}(y)$, $\forall y\in B(0,r)$.
Assume there exists $p_D>0$ such that
\begin{equation}
\label{hypboundedproba}
\forall n>2\rho r\,,\;\forall x\in V_n\,,\quad
\mu_n(\eta_{B(x,r)}\equiv \zeta_{D,x}) \geq p_D\;.
\end{equation}
\item
Assume there exists a function $\psi$ from $\NN\times\NN$ into $\RR$, and a
function $\epsilon$
from $\NN$ to $\RR$, decreasing to $0$, such that if
$B$, $C$ are disjoint subsets of $V_n$, and $e_B$, $e_C$
local events relative to $B$ and $C$ respectively:
\begin{equation}
\label{hypmixing}
|\mu_n(e_B\cap e_C)-\mu_n(e_B)\mu_n(e_C)|
\leq \psi(|B|,|C|)\,\epsilon(dist(B,C))\;.
\end{equation}
\end{enumerate}
Then for any first order sentence $A$,
$$
\lim_{n\rightarrow \infty} \mu_n(A) = 0\mbox{ or } 1\;.
$$
\end{theorem}
Let us notice that the probabilities of local configurations in
(\ref{hypboundedproba}) need not be equal:
stationarity (spatial homogeneity) is not requested. Observe also that
the mixing hypothesis (\ref{hypmixing}) is implied by any strong
mixing condition such as introduced in \cite{Doukhan}.
After the proof of Theorem \ref{th:zuweak} several
examples of random fields satisfying its hypotheses will be given.
\begin{proof}
According to Propositions \ref{prop:reduction} and
\ref{prop:reduction2}, it enough to prove that
the probability of any pattern sentence tends to $1$.
$$
\lim_{n\rightarrow\infty}
\mu_n(\exists x\,D(x))=1\;,
$$
Where $D(0)$ denotes a complete description of the ball
$B(0,r)$. Using the notations of the previous section, let
$\zeta_{D,0}$ be the unique configuration on $B(0,r)$ described by
$D(0)$, and $\zeta_{D,x}$ be its translate on the ball $B(x,r)$:
$$
\forall y\,,\;dist(0,y)\leq r\;,\quad \zeta_{D,x}(x+y)=\zeta_{D,0}(y)\;.
$$
We need to compute the probability that a random configuration
with distribution $\mu_n$ coincides with $\zeta_{D,x}$ on $B(x,r)$,
for at least one $x$. We shall first fix a regular sub-lattice of
$V_n$. Let $R>\rho r$ be a given integer, to be specified later.
Consider the following set of vertices:
\begin{equation}
\label{deftiling}
T_n =
\left\{\,\alpha(2R+1)\,,\;\alpha=0,\ldots,
 \left\lfloor\frac{n}{2R+1}\right\rfloor-1\,\right\}^d\;,
\end{equation}
where $\lfloor\,\cdot\,\rfloor$ denotes the integer part.
Denote by $\tau_n$ the cardinality of $T_n$:
$$
\tau_n = \left\lfloor\frac{n}{2R+1}\right\rfloor^d\;.
$$
Rename the vertices in $T_n$ from $1$ to $\tau_n$:
$x_1,\ldots,x_{\tau_n}$. Since $R>\rho r$, the balls $B(x_i,r)$ and
$B(x_j,r)$ are disjoint for all $i\neq j$. Denote by $e_i$
the event ``the configuration coincides with $\zeta_{D,x_i}$ on
$B(x_i,r)$'' and by $\overline{e_i}=\Xi_n\setminus e_i$ the opposite
event. Obviously, the union from $1$ to $\tau_n$ of the
$e_i$'s implies $\exists x\,D(x)$. We will prove that the
probability of the opposite event tends to $0$.
$$
\lim_{n\rightarrow\infty}
\mu_n\left(\bigcap_{1\leq i\leq \tau_n} \overline{e_i}\right) = 0\;.
$$
For $i=1,\ldots,\tau_n$ denote by $c_i$ the following covariance.
$$
c_i = \mu_n\left(\bigcap_{j=1}^i \overline{e_j}\right)-
\mu_n\left(\bigcap_{j=1}^{i-1} \overline{e_j}\right)
\mu_n\left(\overline{e_i}\right)\;.
$$
An immediate induction leads to:
$$
\mu_n\left(\bigcap_{i=1}^{\tau_n} \overline{e_i}\right) =
c_{\tau_n} +\sum_{i=1}^{\tau_n-1} c_{n-i}\prod_{l=\tau_n-i}^{\tau_n}
\mu_n(\overline{e_l})
+ \prod_{l=1}^{\tau_n}\mu_n(\overline{e_l})\;.
$$
Therefore:
$$
\mu_n\left(\bigcap_{i=1}^{\tau_n} \overline{e_i}\right) \leq
\sum_{i=2}^{\tau_n} |c_i|
+ \prod_{i=1}^{\tau_n}\mu_n(\overline{e_i})\;.
$$
By the first hypothesis of the theorem, the product can be bounded as
follows.
$$
\prod_{l=1}^{\tau_n}\mu_n(\overline{e_l})
\leq (1-p_D)^{\tau_n}\;.
$$
 In order to make sure that it tends to zero, it is enough to
require that $\tau_n$ should tend to infinity. Let us now turn to
the sum of covariances, for which the second hypothesis has to be
used. Observe that by construction of $T_n$, the distance between
any two balls $B(x_i,r)$ and $B(x_j,r)$ is larger than $2(R-\rho
r)/\rho$. Thus the second hypothesis implies that for
$i=1,\ldots,\tau_n$,
$$
\sum_{i=2}^{\tau_n}|c_i|\leq
\sum_{i=2}^{\tau_n}\psi(i\beta(r),\beta(r))\,\epsilon(2(R-\rho r)/\rho)\;.
$$
Let us choose a be a positive increasing real function $\Psi$, tending to
infinity, such that for all $m\in\NN$:
$$
\sum_{i=2}^{m}\psi(i\beta(r),\beta(r))\leq \Psi(m)\;.
$$
Denote by $\Psi^{-1}$ its reciprocal function, which is also positive,
increasing, and tends to infinity.
Define $R=R_n$ as follows.
$$
R_n = \frac{n}{2((\Psi^{-1}(\epsilon(\sqrt{n})^{-1/2}))^{1/d}-1)}\;.
$$
One has~:
$$
\tau_n = \left\lfloor \frac{n}{2R_n+1}\right\rfloor^d
=
\left\lfloor
\frac{1}{\frac{1}{n}+\frac{1}{\Psi^{-1}(\epsilon(\sqrt{n})^{-1/2}))^{1/d}-1}}
\right\rfloor^d\;,
$$
which tends to $+\infty$, since $\epsilon$ tends to zero, and $\Psi^{-1}$ to
$+\infty$.
Now~:
$$
\Psi(\tau_n)\leq \Psi\left(\left(\frac{n}{2R_n}+1\right)^d\right)
=\Psi(\Psi^{-1}(\epsilon(\sqrt{n})^{-1/2}))=\epsilon(\sqrt{n})^{-1/2}\;.
$$
Replacing $\Psi$ by a larger function if necessary,
one can make sure that $2(R_n-\rho r)/\rho\geq \sqrt{n}$.
Since $\epsilon$ is decreasing,
$$
\epsilon(2(R_n-\rho r)/\rho)\leq \epsilon(\sqrt{n})\;,
$$
hence:
$$
\Psi(\tau_n)\,\epsilon(2(R_n-\rho r)/\rho)\leq \epsilon(\sqrt{n})^{1/2}\;,
$$
which tends to $0$ as desired.
\end{proof}
Classical examples of random fields are stationary 
(i.e. translation invariant) random
fields, defined on $\ZZ^d$ as  
sets of random variables $\{\eta(x)\,,\;x\in\ZZ^d\}$ 
(see \cite{Rosenblatt} or \cite{Doukhan}). 
In what follows, we shall view our measure $\mu_n$ as the distribution
of the restriction $\{\eta(x)\,,\;x\in V_n\}$.  
Observe that if the distribution of a stationary random field
$\{\eta(x)\,,\;x\in\ZZ^d\}$ satisfies the
hypotheses of Theorem \ref{th:zuweak}, for $B,C\subset \ZZ^d$, then
the same holds for $B,C\in V_n$ because distances in
the torus $V_n$ are always smaller than in $\ZZ^d$. Notice also, that
even if $\{\eta(x)\,,\;x\in\ZZ^d\}$ is stationary, its restriction to
$V_n$ is not stationary in the sense of 
translations on the torus.

\begin{itemize}
\item {\bf Gaussian models}\\
They are often used in practice, for instance to model  the sea level. Let
$Y=\{Y(x)\,,\;x\in\ZZ^d\}$ be a stationary Gaussian random field, and
$h:\RR\to\{-1,+1\}$ a discrete valued function.
Let $\mu$ be the distribution of $\{h(Y(x))\,,\;x\in\ZZ^d\}$.
Assume that the spectral density $f_Y$ of $Y$
is bounded away from zero. Then in (\ref{hypmixing}), one can choose
$\psi\equiv1$ as a constant and $\epsilon(s)$ as the distance from
$f_Y$ to the set of trigonometric polynomials with degree $<s$.
\item {\bf Ising models}\\
Another classical example of weakly dependent random field is the Ising model
of Definition \ref{def:ising}. In their famous article
\cite{DobrushinShlosman1} proved
the equivalence of several mixing conditions, including
the exponential decay of correlations. Those conditions are satisfied
 in the high temperature case (i.e. for $\mid\!b\!\mid$
small enough) or in large external field (i.e. for $\mid\!a\!\mid$
sufficiently large). See also \cite{Higuchi} for the $2$-dimensional
case. Theorem \ref{th:zuweak} thus yields the zero-one law for the
Ising model with fixed parameters. In Section \ref{ising}, the
parameters $a$ and $b$ will be allowed to depend on $n$.

\item {\bf Weak dependence and mixing}\\
The following setting is described in \cite{DoukhanLouhichi} and we
shall keep their notations.
Let $B$ be a finite subset of $\ZZ^d$. 
Consider the set $\Lambda^{(1)}$ of Lipschitz functions
$\RR\to\RR$ such that
 $$
\mbox{Lip}f=\sup_{x,y\in\{0,1\}^B}\frac{|f(x)-f(y)|}
{\sum_{i\in B}|x_i-y_i|}\le 1,\quad \|f\|_\infty\le 1\;.
$$
The so-called $(\epsilon,\Lambda^{(1)},\phi)$-weak dependence
condition of \cite{DoukhanLouhichi}, defined
from a decreasing sequence $\epsilon_r\searrow 0$
and a function $\phi:\NN^2\times\RR^{+2})\to\RR^+$,
holds if
\begin{equation}\label{depf}
\left|\mbox{Cov}(f_1((X_i)_{i\in B_1}),f_2((X_i)_{i\in B_2})\right|\le
\phi(|B_1|,|B_2|, \mbox{Lip}f_1,\mbox{Lip}f_2)\,\epsilon_{r}
\end{equation}
if $f_1,f_2\in\Lambda^{(1)}$ are as above with $B_1,B_2\subset \ZZ^d$
finite and distant by at least $r$.
A random field satisfying relation (\ref{depf})
for the function $\phi(u,v,U,V)=uU+vV$ 
is called an $\eta$-weakly dependent random field.

In our setting, the $X_t$'s take their values in $\{-1,+1\}$.
In this case, indicators of local configurations may be
written as functions $f((X_i)_{i\in B})$ for a bounded set
$B\subset \ZZ^d$.
Indeed, such functions $f:\{-1,1\}^B\to \{0,1\}$ are
in $\Lambda^{(1)}$ and (\ref{hypmixing}) holds with
$\psi(u,v)=\phi(u,v,1,1)$ and $\epsilon_r$ is as in equation (\ref{depf}).

So in the binary case the weak dependence condition of
\cite{DoukhanLouhichi}
implies the mixing condition (\ref{hypmixing}). In more general
settings, \cite{DoukhanLouhichi} provides a weaker 
alternative to mixing conditions.

Here are some examples of $(\epsilon,\Lambda^{(1)},\phi)$-weakly
dependent fields.

{\it Bernoulli shifts} are defined through an independent and
identically distributed random field $(\xi_t)_{t\in\ZZ^d}$ by
the relation $X_n=H((\xi_{n-t})_{t\in\ZZ^d})$.
If such expressions are (a.s.) well defined,
they generate stationary random fields. If
$$
\delta_r=\EE \left|H((\xi_{t})_{t\in\ZZ^d})-H((\xi_{t}
{\bf 1}_{(|t|\le r)})_{t\in\ZZ^d})\right|
$$
where $|\cdot|$ denotes a norm on $\ZZ^d$,
then it is easy to prove with Doukhan and Louhichi
ideas in \cite{DoukhanLouhichi} that (\ref{depf}) holds with
$\phi(u,v,U,V)=uU+vV$
and $\epsilon(r)=2\delta_{r/2}$. In the case where $H$ takes values in
$\{-1,+1\}$, one has $\psi(u,v)=u+v$.

Other models have similar expressions, in particular
mappings $X_t=h(Y_t)$ of general Bernoulli shifts $(Y_t)$
through functions with values in $\{-1,+1\}.$
As a  simple example, consider a
linear random field, defined by:
$$
Y_t=\sum_{s\in\ZZ^d}a_s\xi_{t-s},
$$
where $\delta_r\le\EE|\xi_0|\sum_{|i|\ge s}|a_i|$.
More general polynomials may also be considered,
they are usually named as Chaotic Volterra expansions
(see \cite{DoukhanLang}).

The next example provides a   way
to exhibit such expansions for $\RR^k$-valued random fields.
It comes from the finance literature where heteroskedasdicity
is a natural feature.
Let  $( \xi_t )_{t \in \ZZ^d}$ be a stationary sequence of random $k \times
m$ matrices, $( a_j )_{j \in \ZZ^{\ast}}$ be a sequence of $m
\times k$ matrices, and $a$ be a vector in $\RR^m$.
A vector valued
$ARCH(\infty)$ random field model is a solution of the recurrence equation
\begin{equation}\label{archvd}
  X_t = \xi_t \left( a + \sum_{j \ne0} a_j X_{t -
  j} \right), \ t\in\ZZ^d
\end{equation}
Such $ARCH(\infty)$ models include a large variety of models (see
\cite{DoukhanTeyssiere}).
 Assume that $\displaystyle\sum_{j \ne0} \| a_j \|\|\xi_0\|_\infty < 1$,
then a stationary  solution
  of equation (\ref{archvd}) in  $L^p$ is given by~:
  \begin{equation}\label{archv2d}
  X_t = \xi_t \left( a + \underset{k = 1}{\overset{\infty}{\sum}}
  \underset{j_1, \ldots, j_k \ne0}{\sum} a_{j_1}\xi_{t
  - j_1} a_{j_2}\ldots a_{j_k} \xi_{t
  - j_1 - \cdots - j_k} a\right) \quad\left(\in\RR^k\right)\;.
\end{equation}
Denote now $A = \underset{j \ne0}{\sum} \| a_j \|$,
$A ( x ) = \underset{\|j \|\geqslant x}{\sum} \| a_j \|$ and $\lambda = \left(
\underset{j \ne0}{\overset{\infty}{\sum}} \| a_j \| \right)  \| \xi_0
\|_\infty$, where 
$$\|(j_1,\ldots,j_k)\|=|j_1|+\cdots+|j_k|\;.$$

In what follows, we assume for simplicity that the random field
$(\xi_t)_{t\in\ZZ^d}$ is iid.

One approximates here $X_t$ by a random variable  independent of $X_0$.
Set $$\overset{}{\tilde{X}_t =} \xi_t \left( a +
\underset{}{\overset{\infty}{\underset{k = 1}{\sum}}}
\underset{\|j_1\| +
\cdots +
\|j_k\| < t}{\sum} a_{j_1} \xi_{t - j_1} \cdots a_{j_k} \xi_{t - j_1 - \cdots -
j_k} a\right) $$
A bound for the error is given by :
  \[ \EE \| X_t - \tilde{X}_t \| \leqslant \EE \| \xi_0 \| \left( \EE \|
     \xi_0 \| \underset{k = 1}{\overset{t - 1}{\sum}} k \lambda^{k - 1} A \left(
     \frac{t}{k} \right) + \frac{\lambda^t}{1 - \lambda} \right) \| a \| \]
 For suitable constants $K$ and  $K'$ we derive
  \[ \| X_t - \tilde{X}_t \| \leqslant \left\{\begin{array}{l}
       K \frac{( \log ( t ) )^{\max\{b , 1\}}}{t^b}  \text{\ if \ }A ( x )
\leqslant C  x^{- b},\
       b, C > 0\\
       K' ( \max\{q, \lambda\} )^{\sqrt{t}}  \text{\ if \ } A ( x ) \leqslant {Cq}^x
     \end{array}\right. \]
Weak dependence of the corresponding random fields follows.
The solution (\ref{archv2d}) of equation (\ref{archvd})
is $\eta-$weakly dependent with
$$\epsilon_r =\| \xi_0 \|_\infty   \| a \|x_r,\mbox{ with }
x_t=  \| \xi_0 \|_\infty
  \sum_{k = 1}^{t - 1} k \lambda^{k - 1} A \left( \frac{t}{k} \right)
  + \frac{\lambda^t}{1 - \lambda} .$$
More precisely, for Riemannian decays of  $\| a_j \|$, with $A ( x )
  \leqslant C  x^{- b}$ for $b, C > 0$, we have
\[
  x_t \leqslant \left\{\begin{array}{l}
        \EE \| \xi_0 \| \left( 2 \EE \| \xi_0 \|
       C \left( \frac{- b}{\ln \lambda} \right)^b \left( \frac{\ln ( t )}{t}
       \right)^b + \frac{\lambda^t}{1 - \lambda} \right) \| a \|
\mbox{\ if \ } b >
       1\\
    \EE \| \xi_0 \| \left( 3 A  C
       \left( \frac{- b}{\ln \lambda} \right) \frac{\ln ( t )}{t^b} +
       \frac{\lambda^t}{1 - \lambda} \right) \| a \| \mbox{\ if \ } 0 < b < 1
     \end{array}\right. \]
For geometric decays $\| a_j \| \le q^j$,
 \[ x_t \leqslant \EE \| \xi_0 \| \left( \underset{k =
     1}{\overset{t - 1}{\sum}} \EE \| \xi_0 \|^k \underset{p =
     t}{\overset{\infty}{\sum}}_{} C^k_p q^p + \frac{\lambda^t}{1 - \lambda}
     \right) \| a \| \]
In those cases, there exist positive constants $K$, $b$, and $C$ 
such that
  \begin{equation*}
    x_t\leqslant \left\{\begin{array}{l}
      K \frac{( \log ( t ) )^{\max\{b,1\}}}{t^b}\,,
\text{ under Riemannian decay  }A ( x ) \leqslant C  x^{-
      b}\\
      K (\max\{ q , \lambda\} )^{\sqrt{t}}\,,
\text{ under geometric decay  } A ( x ) \leqslant Cq^x
    \end{array}\right.
  \end{equation*}

  Now any function $h(X_t)$ of such a random field through a
mapping with values in $\{-1,+1\}$
  satisfies the assumptions of Theorem \ref{th:zuweak}.

Let us finally mention the case of an {\it associated random field},
i.e. a $\{-1,+1\}$-valued random field $\eta$ on $V_n$ such that 
for each coordinatewise non-decreasing function
$f,g:\RR^{V_n}\to\RR$ such that $\EE (f^2(X)+g^2(X))<\infty$, the
covariance between $f(\eta)$ and $g(\eta)$ is non negative
(this is also known as the FKG condition). Such a field can be treated
as before with $\psi(u,v)=uv$ and $\epsilon(s)=\sup_{|u-v|\ge s}\mbox{Cov}(X_u,X_v).$
\end{itemize}
\section{Ising model}
\label{ising}

Here we consider the particular case where the measure $\mu_n$ on the
set of configurations $\Xi_n$ is the Ising model of Definition
\ref{def:ising}:
$$
\forall \eta\in \Xi_n\;,\quad
\mu_{n} (\eta) = \frac{1}{Z}
\exp\left(a(n)\sum_{x\in V_n} \eta(x)
+b(n)\sum_{\{x,y\}\in E_n} \eta(x)\eta(y)\right)\;.
$$
Expectations relative to $\mu_n$ will be denoted by $\EE_n$.
 The surface potential $a=a(n)$ and the pair potential
$b=b(n)$ both depend on $n$. We will assume that $a(n)$ is negative:
we already observed that the model remains unchanged through swapping
$-1$ and $+1$ and changing $a$ to $-a$. If the pair potential $b$
remains bounded, then conditions on $a(n)$ can be given for
the zero-one law.
\begin{theorem}
\label{th:zuising}
Assume that for all $n$, $a(n)<0$ and there exists $b_0>0$
such that $|b(n)|<b_0$. If
$$
 \forall k=1,2,\ldots\;,\quad \lim_{n\rightarrow\infty}
e^{a(n)} n^{\frac{d}{2k}} \in \{ 0\,,\, +\infty\}\;,
$$
then for any first order sentence $A$
$$
\lim_{n\rightarrow\infty}\mu_{n}(A_n) = 0\mbox{ or } 1\;.
$$
\end{theorem}
As a particular case, for $b(n)\equiv 0$, the spins are
independent, $+1$ with probability
$p(n) = e^{a(n)}/(e^{a(n)}+e^{-a(n)})\sim e^{2a(n)}$ or $-1$
with probability $1-p(n)$. The threshold functions for $p(n)$ are
$n^{-\frac{d}{k}}\,,\;k=1,2\ldots$, as expected.
By Proposition \ref{prop:reduction}, it would suffice to prove that
the functions $n^{-\frac{d}{2k}}$ are the threshold functions (for
$e^{a(n)}$) of pattern sentences. Actually, it is possible to
characterize the threshold functions of all basic local sentences.
We begin with the particular case of single pattern sentences.
\begin{proposition}
\label{prop:thresholdpattern}
Assume that for all $n$, $a(n)<0$ and there exists $b_{0}>0$ such that
$|b(n)|<b_{0}$. Let $r$ and $k$ be two integers such that $0 \leq k
\leq \beta(r)$. Let $D(0)$ be a complete description of the ball
$B(0,r)$, with exactly $k$ pluses.
\begin{equation}
\label{threshold1}
\mbox{If } \lim_{n\rightarrow\infty} e^{2a(n)k}n^{d} = 0 \mbox{ then }
\lim_{n\rightarrow\infty} \mu_{n}(\exists x\;D(x)) = 0\;.
\end{equation}
\begin{equation}
\label{threshold2}
\mbox{If } \lim_{n\rightarrow\infty} e^{2a(n)k}n^{d} = +\infty \mbox{
  then }
\lim_{n\rightarrow\infty} \mu_{n}(\exists x\;D(x)) = 1\;.
\end{equation}
\end{proposition}
If $a(n)$ is such that $e^{a(n)}\sim cn^{-\frac{d}{2k}}$, 
\cite{Coupier05} proves that a Poisson approximation holds for
``$\exists x\;D(x)$'' if $D(0)$ has $k$ pluses. Thus
$\mu_n(\exists x\;D(x))$ converges to a limit which is 
different from $0$ and $1$.
\begin{proof}
Let us start with some notations and definitions. Given
$\eta\!\in\!E_{n}=\{-1,+1\}^{V_{n}}$ and $U\!\subset\!V_{n}$,
we denote by $\eta_{U}$ the natural projection over
$\{-1,+1\}^{U}$. If $U$ and $V$ are two disjoint subsets of $V_{n}$
then $\eta_{U}\zeta_{V}$ is the configuration on $U \cup V$ which is
equal to $\eta$ on $U$ and $\zeta$ on $V$. We denote by $\delta U$
the neighborhood of $U$ (corresponding to (\ref{def:voisinage})):
$$
\delta U = \{ y \in V_{n} \setminus U , \; \exists x \in U ,
\; \{ x,y \} \in A_{n} \} ~,
$$
and by $\overline{U}$ the union of the two disjoint sets $U$ and
$\delta U$. In the case of balls, $\overline{B(x,r)} = B(x,r+1)$.
Finally, $\mid\!U\!\mid$ denotes the cardinality of $U$ and
$\mathcal{F}(U)$ the $\sigma$-algebra generated by the configurations
of $\{-1,+1\}^{U}$.\\
Let us fix a complete description $D(0)$ of the ball $B(0,r)$, and
denote by $\zeta_{D,x}$ the local configuration on $B(x,r)$ defined by
the translation of $D(0)$ (notations of Theorem \ref{th:zuweak}).
For any vertex $x$, we define the indicator function $I_{x}^{D}$ as follows:
$$
\forall \eta \in \Xi_{n} , \; I_{x}^{D}(\eta) = \left\lbrace \begin{array}{l}
1 , \; \mbox{ if } \; \eta_{B(x,r)} = \zeta_{D,x} ~,\\
0 , \; \mbox{ otherwise }.
\end{array} \right.
$$
The random variables $I_x^D$ indicate occurrences of a given local
configuration on the lattice. Their
sum counts these occurrences. It will be denoted by $X_n^D$.
$$
X_{n}^{D} = \sum_{x \in V_{n}} \; I_{x}^{D}\;.
$$
Due to periodicity, this sum bears over $|V_{n}|=n^{d}$
terms. The event $X_{n}^{D} > 0$ is
logically equivalent to the single pattern sentence $\exists
x\;D(x)$. So, the quantity we are actually interested in is
$\mu_{n} (X_{n}^D > 0)$.

Let us now fix a vertex $x$ and simply denote by $B$ the ball
$B(x,r)$. Recall that there exists a finite number of configurations
on $B$ (precisely $2^{\beta(r)}$).
Denote by $\mathcal{D}_{r}^{x}$ their set.
The \textit{local energy} $H^{B}(\sigma)$ of the
configuration $\sigma \in \{-1,+1\}^{\overline{B}}$
on  $\overline{B}$ is defined by:
\begin{equation}
\label{local_energy}
H^{B}(\sigma) = a(n) \sum_{y \in B} \sigma(y) + b(n)
\sum_{{\scriptstyle \lbrace y,z \rbrace \in A_{n}} \atop
\scriptstyle (y \in B) \vee (z \in B)} \sigma(y) \sigma(z) ~.
\end{equation}
Hence, for all $\sigma \in \{-1,+1\}^{\delta B}$, the local energy
$H^{B}(\zeta_{D,x} \sigma)$ on $B$ of the configuration which is equal
to
$\zeta_{D,x}$ on $B$ and $\sigma$ on $\delta\!B$, can be expressed as:
$$
H^{B}(\zeta_{D,x} \sigma) = a(n) (2k - \beta(r)) + b
\left( \sum_{{\scriptstyle \{y,z\} \in A_{n}} \atop \scriptstyle y,z
    \in B}
\zeta_{D,x}(y) \zeta_{D,x}(z) + \sum_{{\scriptstyle \{y,z\} \in A_{n}}
  \atop
\scriptstyle y \in B, z \in \delta B} \zeta_{D,x}(y) \sigma(z) \right) ~,
$$
where $k$ is the number of $+1$'s of the complete description
$D(0)$.
This notion of local energy naturally appears in the expression of the
conditional probability $\mu_{n}(I_{x}^{D}=1\!\mid\!\sigma)$,
$\sigma \in \{-1,+1\}^{\delta\!B}$:
\begin{equation}
\label{proba_cond}
\mu_{n}(I_{x}^{D}=1 \mid \sigma) = \frac{e^{H^{B}(\zeta_{D,x}\sigma)}}
{\sum_{\zeta' \in \mathcal{D}_{r}^{x}} e^{H^{B}(\zeta' \sigma)}} ~.
\end{equation}
Moreover, it is connected to the number $X_{n}^{D}$ of copies of
$\zeta_{D,0}$ in the graph since:
$$
\EE_{n} \lbrack X_{n}^{D} \rbrack = \EE_{n} \left\lbrack n^{d}
\mu_{n}(I_{x}^{D}=1 \mid \mathcal{F}(\delta B)) \right\rbrack ~.
$$
Here, $\mu_{n}(I_{x}^{D}=1\!\mid\!\mathcal{F}(\delta\!B))$
represents a $\mathcal{F}(\delta\!B)$-measurable random variable,
whereas for $\sigma \in \{-1,+1\}^{\delta\!B}$,
$\mu_{n}(I_{x}^{D}=1\!\mid\!\mathcal{F}(\delta B))(\sigma) =
\mu_{n}(I_{x}^{D}=1\!\mid\!\sigma)$ is a conditional probability.
\vskip 3mm
Note that the set $\delta\!B$ has a bounded cardinality.
Thus, the pair potential $b(n)$ being bounded, there exist two
constants $m$ and $M$ (which do not depend on $n$) such that for all
local configuration $\zeta' \in \mathcal{D}_{r}^{x}$ having exactly
$k'$ pluses and for all configuration $\sigma \in \Xi_{n}$:
\begin{equation}
\label{encadrement}
2a(n) (k' - k) + m \leq H^{B}(\zeta' \sigma_{\delta B}) -
H^{B}(\zeta_{D,x} \sigma_{\delta B}) \leq 2a(n) (k' - k) + M ~.
\end{equation}
Statements (\ref{threshold1}) and (\ref{threshold2}) rely on the fact
that the difference between the local energies
$H^{B}(\zeta' \sigma_{\delta\!B})$ and $H^{B}(\zeta_{D,x}
\sigma_{\delta\!B})$ is uniformly bounded, for any configuration on
$\delta\!B$.
\vskip 1mm
Assume that $k>0$ and $\lim
e^{2a(n)k}n^{d}=0$. Let us denote by $\zeta^{0}$ the unique element of
$\mathcal{D}_{r}^{x}$ having no $+1$. Relations
(\ref{proba_cond}) and (\ref{encadrement}) imply that, for all
$\sigma \in \Xi_{n}$,
\begin{eqnarray*}
\mu_{n}(I_{x}^{D}=1 \mid \sigma_{\delta B}) & \leq &
e^{H^{B}(\zeta_{D,x}
\sigma_{\delta B}) - H^{B}(\zeta^{0} \sigma_{\delta B})} \\
& \leq & e^{2a(n)k - m} ~.
\end{eqnarray*}
As a consequence:
\begin{eqnarray*}
\EE_{n} \lbrack X_{n}^{D} \rbrack & = & n^{d} \EE_{n} \left\lbrack
\mu_{n}(I_{x}^{D}=1 \mid \mathcal{F}(\delta B)) \right\rbrack \\
& \leq & n^{d} e^{2a(n)k - m} ~,
\end{eqnarray*}
which tends to $0$ as $n$ tends to infinity. Finally, $X_{n}^{D}$
being an integer valued variable, its probability of being positive
is bounded by its expectation and relation (\ref{threshold1}) follows.
\vskip 1mm
Assume now that $\lim e^{2a(n)k}n^{d} =+\infty$. 
For all local configurations $\zeta' \in
\mathcal{D}_{r}^{x}$ and for all $\sigma \in \Xi_{n}$,
(\ref{encadrement}) implies that the difference $H^{B}(\zeta'
\sigma_{\delta\!B})-H^{B}(\zeta_{D,x} \sigma_{\delta\!B})$ is bounded
by $-2a(n)k+M$ (the surface potential $a(n)$ is negative). Hence,
using relation (\ref{proba_cond}), there exists a positive constant
$c=2^{-\beta(r)}e^{-M}$ satisfying:
\begin{equation}
\label{minor_proba_cond}
\forall \sigma \in \Xi_{n}, \; \mu_{n}(I_{x}^{D}=1 \mid \sigma_{\delta
  B})
\geq c e^{2a(n)k} ~.
\end{equation}
As in the proof of Theorem \ref{th:zuweak}, let us introduce the subset
$T_n$ of $V_n$ defined by:
$$
T_n =
\left\{\,\alpha(\rho(2r+1)+1)\,,\;\alpha=0,\ldots,
 \left\lfloor\frac{n}{\rho(2r+1)+1}\right\rfloor-1\,\right\}^d\;,
$$

\begin{figure*}[!ht]
\label{fig:reseau}
\begin{center}
\includegraphics[width=8cm,height=5cm]{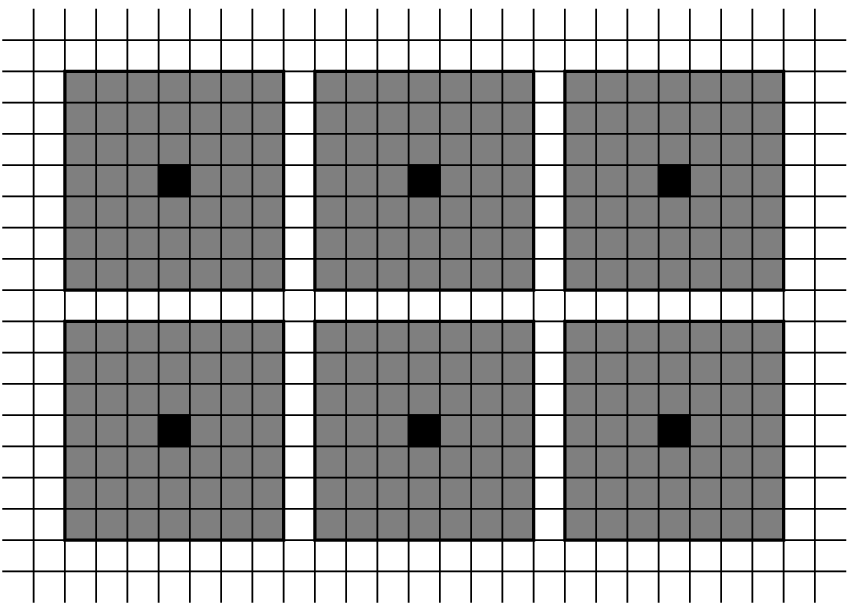}
\caption{The grey vertices represent the set
$\mathcal{T}_{n}$ (in dimension $d=2$, with $r=3$, $\rho=1$,
and $p=+\infty$). The black vertices are
 elements of $T_{n}$.}
\end{center}
\end{figure*}

Denote by $\tau_{n}$ the cardinality of $T_{n}$ and $\mathcal{T}_{n}$
the union of the balls of radius $r$, centered at the elements of
$T_{n}$ (see Figure \ref{fig:reseau}).
Firstly, let us remark that there exists a positive constant $c'$
satisfying $\tau_{n} \geq c' n^{d}$. Secondly, if $x$ and $x'$
are two distinct elements of $T_{n}$ then no vertex of the ball
$B(x,r)$ can be neighbor with a vertex of the ball $B(x',r)$. Thus, let us
define the random variable $\tilde{X}_{n}^{D}$ by:
$$
\tilde{X}_{n}^{D} = \sum_{x \in T_{n}} \; I_{x}^{D} ~.
$$
As $\tilde{X}_{n}^{D} \leq X_{n}^{D}$, it suffices to prove that
$\mu_{n} (\tilde{X}_{n}^{D}=0)$ tends to $0$ (to obtain
(\ref{threshold2})). 
The Gibbs measure $\mu_n$ is a Markov random field with respect to the
graph structure defined by (\ref{def:voisinage}). This has two
consequences. On the one hand, due to the definition of $T_n$, one has:
\begin{eqnarray}
\label{indep_cond}
\mu_{n}\left(\tilde{X}_{n}^{D} = 0 \mid \mathcal{F}
(\delta \mathcal{T}_{n})\right)
& = & \mu_{n}\left( \bigwedge_{x \in T_{n}}
I_{x}^{D}=0 \mid \mathcal{F}(\delta \mathcal{T}_{n}) \right) \nonumber \\
& = & \prod_{x \in T_{n}} \mu_{n}\left( I_{x}^{D}=0 \mid
\mathcal{F}(\delta \mathcal{T}_{n}) \right) ~.
\end{eqnarray}
On the other hand, for all $x\in T_n$, the random variable
$\mu_{n}(I_{x}^{D}=0\!\mid\!\mathcal{F}(\delta \mathcal{T}_{n}))$ is
$\mathcal{F}(\delta\!B(x,r))$-measurable:
\begin{equation}
\label{markov1}
\mu_{n}\left(I_{x}^{D}=0 \mid \mathcal{F}(\delta
  \mathcal{T}_{n})\right) =
\mu_{n} \left(I_{x}^{D}=0 \mid \mathcal{F}(\delta B(x,r))\right) ~.
\end{equation}
Let $\sigma^{max}$ be an element of $\Xi_{n}$ maximizing the conditional
probability $\mu_{n}(I_{x}^{D}=0\!\mid\!\sigma_{\delta\!B(x,r)})$.
Then, from (\ref{indep_cond}) and (\ref{markov1}), we get for any
vertex $x$ and for any $\sigma \in \Xi_{n}$:
\begin{eqnarray*}
\mu_{n}\left(\tilde{X}_{n}^{D} = 0 \mid
\mathcal{F}(\delta \mathcal{T}_{n})\right) & \leq & \left(
\mu_{n}(I_{x}^{D}=0
\mid \sigma_{\delta B(x,r)}^{max}) \right)^{\tau_{n}} \\
& \leq & e^{- \tau_{n} \mu_{n}\left(I_{x}^{D}=1 \mid \sigma_{\delta B(x,r)}^{max}\right)} ~.
\end{eqnarray*}
Finally, the inequality (\ref{minor_proba_cond}) implies that:
\begin{eqnarray*}
\mu_{n}\left(\tilde{X}_{n}^{D} = 0\right) & = & \EE_{n}
\left\lbrack
\mu_{n}\left(\tilde{X}_{n}^{D} = 0 \mid \mathcal{F}(\delta
\mathcal{T}_{n})\right) \right\rbrack \\
& \leq & e^{- \tau_{n} \mu_{n}\left(I_{x}^{D}=1 \mid
\sigma_{\delta B(x,r)}^{max}\right)} \\
& \leq & e^{-c\;c' n^{d} e^{2a(n)k}}
\end{eqnarray*}
which tends to $0$ as $n$ tends to infinity by hypothesis.
\end{proof}

Proposition \ref{prop:thresholdpattern} shows that the appearance of a
given local pattern only depends on its number of $+1$'s:
if $e^{a(n)}$ is small compared to $n^{-\frac{d}{2k}}$, then no
pattern of fixed size, with $k$ pluses, should appear in the
configuration. If $e^{a(n)}$ is large compared to $n^{-\frac{d}{2k}}$, all
configurations with $k$ pluses should appear.
In other words, the
threshold function for the appearance of a given pattern is
$n^{-\frac{d}{2k}}$, where $k$ is the number of $+1$'s in the
pattern. Proposition \ref{prop:basiclocal} below will show that the
threshold function for a basic local sentence $L$ is
$n^{-\frac{d}{2k(L)}}$, where $k(L)$ is an integer that we call the
\textit{index} of $L$. Its definition refers to the decomposition
(\ref{decompositionlocalcomplete}) of a local property into a finite
disjunction of complete descriptions, already used in the proof of
Proposition \ref{prop:reduction}.

\begin{definition}
\label{def:index}
Let $L$ be the basic local sentence defined by:
$$
\exists x_{1}\ldots\exists x_{m}\; \left( \bigwedge_{1\leq i<j\leq m}
dist(x_{i},x_{j})>2r \right)\;\wedge\;\left( \bigwedge_{1\leq i \leq m}
\psi_{i}(x_{i}) \right)\;.
$$
If $L$ is not satisfiable, then we shall set $k(L)=+\infty$. If
$L$ is satisfiable, for each $i=1,\ldots,m$, consider the finite set
$\{ D_{i,1},\ldots,D_{i,d_{i}}\}$ of those complete descriptions on
the ball $B(x_{i},r)$ which imply $\psi_{i}(x_{i})$.
$$
\psi_{i}(x_{i}) \leftrightarrow \bigvee_{1\leq j\leq d_{i}} D_{i,j}(x_{i})\;.
$$
Each complete description $D_{i,j}(x_i)$ corresponds to a
configuration on $B(x_{i},r)$. Denote by $k_{i,j}$ its number of $+1$'s.\\
The {\rm index} of $L$, denoted by $k(L)$ is defined by:
\begin{equation}
\label{defindex}
k(L) = \max_{1\leq i\leq m} \min_{1\leq j\leq d_{i}} k_{i,j} ~.
\end{equation}
\end{definition}
The intuition behind Definition \ref{def:index} is the following.
Assume $e^{a(n)}$ is small compared to $n^{-\frac{d}{2k(L)}}$. Then
there exists $i$ such that none of the $D_{i,j}(x_{i})$ can be
satisfied, therefore there is no $x$ such that $\psi_{i}(x)$
is satisfied, and $L$ cannot be satisfied either. 
On the contrary, if $e^{a(n)}$
is large compared to $n^{-\frac{d}{2k(L)}}$, then for all
$i=1,\ldots,m$, $\psi_{i}(x_{i})$ should be satisfied for at least one
vertex $x_{i}$, and the probability of satisfying $L$ should be
large. In other words, $n^{-\frac{d}{2k(L)}}$ is the threshold function of $L$.

\begin{proposition}
\label{prop:basiclocal}
Assume that for all $n$, $a(n)<0$ and there exists $b_{0}>0$ such that
$|b(n)|<b_{0}$.  
Let $L$ be a basic local property, and $k(L)$ be its index. If $L$ is
satisfiable and $k(L)>0$, then its threshold function for $e^{a(n)}$
is $n^{-\frac{d}{2k(L)}}$. If $k(L)=0$, its probability tends to $1$.
\end{proposition}
\begin{proof}
Assume $L$ is satisfiable (otherwise its probability is null). For
$a(n)<0$, we need to prove that $\mu_{n}(L)$ tends to $0$ if
$e^{2 a(n) k(L)} n^{d}$ tends to $0$ (in this case, $k(L) > 0$), and
that it tends to $1$ if $e^{2 a(n) k(L)} n^{d}$ tends to $+\infty$.
The former will be proved first. Consider again the decomposition of
$L$ into complete descriptions:
$$
L \leftrightarrow \exists x_{1}\ldots\exists x_{m}\;
\left( \bigwedge_{1\leq i<j\leq m} dist(x_{i},x_{j})>2r \right)
\;\wedge\;
\left( \bigwedge_{1\leq i \leq m} \bigvee_{1\leq j\leq d_{i}}
D_{i,j}(x_{i}) \right) ~.
$$
If $e^{2 a(n) k(L)} n^{d}$ tends to $0$, there exists $i$ such that:
$$
\forall j=1,\ldots, d_{i}\,,\;
\lim_{n \to \infty} e^{a(n)}n^{\frac{d}{2k_{i,j}}} = 0\;.
$$
By Proposition \ref{prop:thresholdpattern}, the probability of
$(\exists x\;D_{i,j}(x))$ tends to $0$ for all $j=1,\ldots,d_{i}$.
Therefore the probability of $(\exists x\,\psi_{i}(x))$ tends to $0$,
which implies that $\mu_{n}(L)$ tends to $0$.\\
Conversely, for each $i=1, \ldots ,m$, choose one of the $D_{i,j}$'s,
such that the number of $+1$'s in the corresponding configuration
is minimal (among all $k_{i,j}$'s). Denote that particular complete
description by $D_{i}$ and by $k_{i}$ its number of $+1$'s
(hence, $k(L) = \max_{1 \leq i \leq m} k_{i}$). As in the proof of
Proposition \ref{prop:thresholdpattern}, we shall use the lattice
$T_{n}$. Remember that its cardinality $\tau_{n}$ is of order
$n^{d}$. The basic local sentence $L$ is implied by the sentence
$L'$ defined as follows:
$$
L' \leftrightarrow \exists x_{1} \ldots \exists x_{m} \;
\left( \bigwedge_{1 \leq i \leq m} x_{i} \in T_{n} \right) \;
\wedge \; \left(\bigwedge_{1 \leq i < j \leq m} x_{i} \not= x_{j}
\right)
\; \wedge \; \left( \bigwedge_{1 \leq i \leq m} D_{i}(x_{i}) \right)~,
$$
For $i = 1, \ldots ,m$, let us introduce the random variable
$\tilde{X}_{n}^{D_{i}}$~:
$$
\tilde{X}_{n}^{D_{i}} = \sum_{x \in T_{n}} \; I_{x}^{D_{i}} ~.
$$
Even if some of the $D_{i}$'s correspond to the same complete
description, the event
$$
\bigwedge_{i=1}^{m} \;\left\lbrace \tilde{X}_{n}^{D_{i}} \geq m \right\rbrace
$$
implies the sentence $L'$ (and also the basic local sentence $L$).
As a consequence, it is enough to prove the convergence of
$\mu_{n}(\tilde{X}_{n}^{D_{i}} \geq m)$ to $1$ for any index
$i = 1, \ldots ,m$.\\
Fix such an index $i$. We have already seen in the proof of Proposition
\ref{prop:thresholdpattern} that, for two different vertices $x$ and
$x'$ of $T_{n}$,
$$
\begin{array}{ll}
(i) \quad \mu_{n}\left(I_{x}^{D_{i}}=1 \mid
\mathcal{F}(\delta \mathcal{T}_{n})\right) =
\mu_{n}\left(I_{x'}^{D_{i}}=1
\mid \mathcal{F}(\delta \mathcal{T}_{n})\right)~; \\
(ii) \quad \mu_{n}\left(I_{x}^{D_{i}}=I_{x'}^{D_{i}}=1
\mid \mathcal{F}(\delta \mathcal{T}_{n})\right) =
\mu_{n}\left(I_{x}^{D_{i}}=1 \mid \mathcal{F}(\delta
\mathcal{T}_{n})\right) \; \times \; \mu_{n}
\left(I_{x'}^{D_{i}}=1 \mid \mathcal{F}(\delta \mathcal{T}_{n})\right)~.
\end{array}
$$
These two statements mean that, conditioned on a given configuration
$\sigma$ on $\delta\!\mathcal{T}_{n}$, the distribution of the random variable
$\tilde{X}_{n}^{D_{i}}$ is binomial with
parameters $\tau_{n}$ and $\mu_{n}(I_{x}^{D_{i}}=1\!\mid\!\sigma)$.
Moreover, for any vertex $x$, the inequality (\ref{minor_proba_cond}) says:
\begin{equation}
\label{uniformconfig}
\exists \; c^{''} > 0 , \; \forall \; \sigma \in \Xi_{n}, \; \tau_{n} \;
\mu_{n}\left(I_{x}^{D_{i}}=1 \mid \sigma_{\delta
    \mathcal{T}_{n}}\right)
\geq c^{''} n^{d} e^{2 a(n) k_{i}} ~.
\end{equation}
Now, assume that $n^{d} e^{2 a(n) k(L)}$ tends to infinity. Since
$k_{i} \leq k(L)$ and $a(n)$ is negative, relation
(\ref{uniformconfig})
involves the convergence of the product $\tau_{n} \times
\mu_{n}\left(I_{x}^{D_{i}}=1\!\mid\!\sigma_{\delta\!\mathcal{T}_{n}}\right)$
to infinity, uniformly on $\sigma \in \Xi_{n}$. Then, it is easy to
check that the conditional probability
$\mu_{n}\left(\tilde{X}_{n}^{D_{i}}
\geq m\!\mid\!\sigma_{\delta\!\mathcal{T}_{n}}\right)$ tends to $1$,
uniformly on the configuration $\sigma \in \Xi_{n}$. Therefore the probability
$$
\mu_{n}\left(\tilde{X}_{n}^{D_{i}} \geq m \right) =
\EE_{n} \left\lbrack \mu_{n}\left(\tilde{X}_{n}^{D_{i}}
\geq m \mid \mathcal{F}(\delta \mathcal{T}_{n})\right) \right\rbrack
$$
tends to $1$ as $n$ tends to infinity.
\end{proof}

Having characterized the threshold functions of all basic local
properties, the proof of Theorem \ref{th:zuising} is now clear. If
$e^{a(n)}n^{\frac{d}{2k}}$ tends to $0$ or $+\infty$ for any positive
integer $k$, then by Proposition \ref{prop:basiclocal} the probability
of any basic local sentence tends to $0$ or $1$. This remains true for
any boolean combination of basic local sentences (cf. Proposition
\ref{prop:reduction}). By Gaifman's theorem, these boolean
combinations
cover all first-order sentences. Hence the zero-one law holds for
 first-order logic.

\bibliographystyle{alea2}
\bibliography{zero_un}
\end{document}